\documentclass[twoside,12pt,leqno]{article}
\advance\topmargin by -\headheight\advance\topmargin by -\headsep
\newcommand{\tmvolxx}{xx}
\newcommand{\tmyearyyyy}{yyyy}

\newcommand{\FirstPageHead}[3]{
{\footnotesize
\vskip -8mm
\centerline {Travaux math\'ematiques, \quad
Volume #1 (#2),
#3,\quad \copyright\  Universit\'e du Luxembourg}}\vspace{-3mm}}

\usepackage{amsmath}
\usepackage{amsthm}
\usepackage{amssymb}
\usepackage{amscd}
\usepackage{graphicx}
\usepackage{epsfig}
\usepackage{amssymb, amsmath, amsfonts, amscd}
\input xypic
\numberwithin{equation}{section}
\newtheorem{theorem}{Theorem}[section]

\newtheorem{proposition}[theorem]{Proposition}
\newtheorem{corollary}[theorem]{Corollary}

\theoremstyle{definition}
\newtheorem{definition}[theorem]{Definition}

\numberwithin{equation}{section}

\allowdisplaybreaks
\begin{document}
\thispagestyle{empty}
\FirstPageHead{\tmvolxx}{\tmyearyyyy}{\pageref{firstpage}--\pageref{lastpage}}
\label{firstpage}

\newcommand{\m}{\mathfrak{m}}
\newcommand{\F}{\mathcal{F}}
\newcommand{\G}{\mathcal{G}}
\newcommand{\I}{\mathcal{I}}
\newcommand{\A}{\mathcal{A}}
\newcommand{\B}{\mathcal{B}}
\renewcommand{\L}{\mathcal{L}}
\newcommand{\E}{\mathcal{E}}
\renewcommand{\O}{\mathcal{O} }
\renewcommand{\a}{\alpha}
\renewcommand{\b}{\beta }
\newcommand{\g}{\gamma}
\renewcommand{\P}{\mathbf{P} }
\renewcommand{\Pr}{\mathcal{P} }

\renewcommand{\u}{\frac{1}{u} }
\renewcommand{\t}{\frac{1}{t} }
\newcommand{\mf}[1]{\mathfrak{#1} }
\newcommand{\mc}[1]{\mathcal{#1} }
\newcommand{\vw}{\mathbf{V}_W}
\newcommand{\pa}{\partial}
\newcommand{\onab}{\overline{\nabla}}

\renewcommand{\equiv}{\backsim }
\renewcommand{\o}[1]{\overline{#1}}
\newcommand{\Aut}{\operatorname{Aut}} 
\newcommand{\SL}{\operatorname{SL}} 
\newcommand{\GL}{\operatorname{GL}} 
\newcommand{\Der}{\operatorname{Der} }
\newcommand{\diff}{\operatorname{Diff}} 
\newcommand{\Hom}{\operatorname{Hom} }
\newcommand{\Ext}{\operatorname{Ext} }
\newcommand{\End}{\operatorname{End} }
\newcommand{\Spec}{\operatorname{Spec} }
\newcommand{\Tor}{\operatorname{Tor} }
\renewcommand{\H}{\operatorname{H} }
\newcommand{\Proj}{\operatorname{Proj} }
\newcommand{\D}{\operatorname{D} }
\newcommand{\HH}{\operatorname{HH} }
\renewcommand{\S}{\operatorname{Sym} }
\newcommand{\K}{\operatorname{K} }
\newcommand{\C}{\operatorname{C} }
\newcommand{\R}{\operatorname{R} }
\renewcommand{\K}{\operatorname{K} }
\renewcommand{\lg}{\mathfrak{g}}
\newcommand{\bn}{B^{\otimes n} }
\newcommand{\s}{\sigma}
\newcommand{\Z}{\mathbf{Z} }
\newcommand{\sym}{\operatorname{Sym} }
\newcommand{\e}{\overline{e}}
\newcommand{\Diff}{\operatorname{Diff}}
\newcommand{\p}{\mathbf{P}}
\renewcommand{\mf}[1]{\mathfrak{#1} }
\renewcommand{\mc}[1]{\mathcal{#1} }
\newcommand{\mb}[1]{\mathbf{#1} }




\markboth{Helge Maakestad}{Jet bundles on projective space}
$ $
\bigskip

\bigskip

\centerline{{\Large   Jet bundles on projective space}}

\bigskip
\bigskip
\centerline{{\large by  Helge Maakestad}}

\vspace*{.7cm}

\begin{abstract}
In this paper we state and prove some 
  results on the structure of the jetbundles as left and right module
  over the structure sheaf $\O$ on the projective line and projective
  space using elementary techniques involving diagonalization of
  matrices, multilinear algebra and sheaf cohomology.

\end{abstract}

\pagestyle{myheadings}
\section{Introduction}

In  this paper old and new results on the structure of the jetbundles
on the projective line and projective space as left and right module
is given. In the first section explicit techniques and known results on the splitting type
of the jetbundles as left and right module over the structure sheaf
$\O$ are recalled. In the final section some new results on
the structure of the jetbundles on $\p^N=\SL(V)/P$ as left and right $P$-module are
proved using sheaf cohomology and multilinear algebra. 

\section{On the left and right structure and matrix diagonalization}

In this section we recall techniques and results introduced in previous
papers (\cite{MAA}, \cite{maa2}) where the jets are studied as left and right
module using explicit calculations involving diagonalization of matrices.

\textbf{Notation}. Let $X/S$ be a separated scheme, and let $p,q$
be the two projection maps $p,q:X\times X\rightarrow X$. There is a
closed immersion
\[ \Delta :X\rightarrow X\times X \]

and an exact sequence of sheaves 
\begin{equation}  \label{prinexact} 
  0\rightarrow \I^{k+1} \rightarrow \O_{X\times X} \rightarrow \O_{\Delta^k} \rightarrow 0. 
\end{equation} 
The sheaf $\I$ is the sheaf defining the diagonal in $X\times X$.

\begin{definition} \label{principalparts} 
Let $\E$ be a locally free  $\O_X$-module. We define
the \emph{$k$'th order jets} of $\E$, to be
\[ \Pr^k_X(\E)=p_*(\O_{\Delta^k}\otimes_{X\times X} q^*\E).\]
Let $\Pr^k_X$ denote the module $\Pr^k_X(\O_X)$.
\end{definition}

There is the following result:

\begin{proposition} \label{fundamentalexact} Let $X/S$ be smooth and let $\E$ be a locally free
$\O_X$-module.
There exists an exact sequence
\[ 0 \rightarrow S^k(\Omega^1_X)\otimes \E \rightarrow \Pr^k_X(\E) \rightarrow
\Pr^{k-1}_X(\E) \rightarrow 0 \]
of left $\O_X$-modules, where $k=1,2,\dots$.
\end{proposition}
\begin{proof} For a proof see \cite{perk1}, section 4. 
\end{proof}

It follows that for a smooth morphism $X\rightarrow S$ of relative dimension
$n$, and $\E$ a locally free sheaf on $X$ of rank $e$, the jets
$\Pr^k(\E)$ is locally free of rank $e\binom{n+k}{n}$.

\begin{proposition} \label{functoriality}
Let $f:X\rightarrow Y$ be a map of smooth schemes
over $S$, and let 
$\E$ be a locally free $\O_Y$-module. 
There exists a commutative diagram of exact sequences
\[ \diagram 0 \rto & S^k(f^*\Omega^1_Y)\otimes f^*\E \dto \rto
& f^*\Pr^k_Y(\E) \dto \rto & f^*\Pr^{k-1}_Y(\E) \dto \rto & 0 \\
0 \rto & S^k(\Omega^1_X)\otimes f^*\E \rto & \Pr^k_X(f^*\E) \rto & \Pr^{k-1}_X(f^*\E)
\rto & 0 \enddiagram \]
of left $\O_X$-modules for all $k=1,2,\dots$.
\end{proposition}
\begin{proof} See \cite{perk1}. 
\end{proof}
Hence for any open set $U\subseteq X$ of a scheme $X$ there is an isomorphism
\begin{equation}\label{local}
\Pr^k_X(\E)|_U\cong \Pr^k_U(\E|_U)
\end{equation}
of left $\O_U$-modules.

Given vector bundles $\F$ and $\G$ one may define the \emph{sheaf of
  polynomial differential operators} of order $k$
from $\F$ to $\G$ (following \cite{egaIV} section 16.8), denoted
$\Diff_X ^k(\F,\G)$. There is an isomorphism
\begin{equation} \label{prindiff}
  Hom _X(\Pr_X ^k( \F ),\G)\cong \Diff_X ^k(\F ,\G ) 
\end{equation} 
of sheaves of abelian groups. 
Hence the sheaf of jets $\Pr_X^k(\F)$ is in a natural way a sheaf of $\O_X$-bimodules.
We write $\Pr^k(\E)^{L}$ (resp. $\Pr^k(\E)^{R}$) to indicate we are
considering the left (resp. right) structure.
Note that for $X$ smooth over $S$ and $\E$ locally free
of finite rank, it follows that $\Pr^k_X(\E)$ is locally free of finite rank as left and right $\O_X$-module separately.
%
%
%
%
%
%

%
%
%
%
The aim of this section is to give formulas for the transition-matrices $L^k_n$ and $R^k_n$ defining the jets
 $\Pr^k(\O(n))$ as locally free sheaves on the projective line.
The structure-matrices defining the jets on the projective line are matrices with entries 
Laurent-polynomials with binomial coefficients. The matrix $L^k_n$
define the left $\O$-module structure, and the matrix $R^k_n$ define the right $\O$-module structure.

%
%
%
%
%

\textbf{Notation}.
Consider the projective line $\P=\P^1_\Z$ over the integers $\Z$. By definition $\P$ equals $\Proj \Z[x_0,x_1]$. Let 
$U_i=\Spec \Z[x_0,x_1]_{(x_i)}$. Let furthermore
$t=\frac{x_1}{x_0}\otimes 1$ and $u=1\otimes \frac{x_1}{x_0}$ under the isomorphism 
 \[ \Z[\frac{x_1}{x_0}]\otimes_{\Z} \Z[\frac{x_1}{x_0}]\cong \Z[u,t] .\]
It follows from Definition \ref{principalparts}  and Proposition \ref{local} that the sheaf of jets 
 $\Pr^k(\O(n))|_{U_0}$ as a left $\O_{U_0}$-module equals
\[ \Z[t]\{1\otimes x_0^n, dt\otimes x_0^n, \dots , dt^k\otimes x_0^n\} \]
as a free $\Z[t]$-module where $dt^p=(u-t)^p$. 
Let 
$\t=\frac{x_0}{x_1}\otimes 1$ and $\u=1\otimes \frac{x_0}{x_1}$ under the isomorphism 
 \[ \Z[\frac{x_0}{x_1}]\otimes_{\Z} \Z[\frac{x_0}{x_1}]\cong \Z[\u,\t] .\]
It follows from Definition \ref{principalparts}  and Proposition \ref{local} that  
$\Pr^k(\O(n))|_{U_1}$ as a left $\O_{U_1}$-module equals
\[ \Z[s]\{1\otimes x_1^n,ds\otimes x_1^n, \dots , ds^k\otimes x_1^n\} \]
as a free $\Z[s]$-module, where $s=\frac{1}{t}$ and $ds^p=(\u-\t)^p$. 
Let 
\[ C=\{1\otimes x_0^n,dt\otimes x_0^n, \dots , dt^k\otimes x_0^n\} \]
and 
\[ C'=\{1\otimes x_1^n,ds\otimes x_1^n, \dots , ds^k\otimes x_1^n\} .\]
We aim to describe the jets by calculating the transition-matrix 
\[ L^k_n=[I]^{C'}_C\]
between
the bases $C'$ and $C$ as an element of the group $\GL(k+1,\Z[t,t^{-1}])$. 
The following theorem gives the structure-matrix $L^k_n$ describing the left $\O$-module structure
of the sheaf of jets on the projective line over $\Z$. 
\begin{theorem} \label{leftstructure} Consider $\Pr^k(\O(n))$ as left $\O$-module on $\P$. 
In the case $n<0$ the structure-matrix $L^k_n$ is given by the following formula:
\begin{equation}\label{left1}
ds^p\otimes x_1^n=\sum_{j=p}^{k}(-1)^{j}\binom{j-n-1}{p-n-1}t^{n-p-j} dt^{j}\otimes x_0^n ,
\end{equation}
where $0\leq p \leq k$. If $n\geq 0$ the matrix $L^k_n$ is given as follows:

Case $1\leq k \leq n$: 
\begin{equation} \label{left2}
ds^p\otimes x_1^n= \sum_{j=p}^{k}(-1)^p \binom{n-p}{j-p}t^{n-p-j}dt^{j}\otimes x_0^n, 
\end{equation}
where $0\leq p \leq k$. 

Case $n<k$: 
\begin{equation}\label{left3}
ds^p\otimes x_1^n=\sum_{j=p}^{n}(-1)^p\binom{n-p}{j-p}t^{n-p-j}dt^{j}\otimes x_0^n, 
\end{equation}
if $1\leq p \leq n$, and 
\begin{equation}\label{left4}
ds^p\otimes x_1^n=\sum_{j=p}^{k}(-1)^{j}\binom{j-n-1}{p-n-1}t^{n-p-j}dt^{j}\otimes x_0^n, 
\end{equation}
if $n<p\leq k$. 
\end{theorem} 
\begin{proof} See \cite{maa2}.
\end{proof} 
%
%
%
We use similar calculations as in the Theorem \ref{leftstructure} to study the structure-matrix $R^k_n$ defining the jets
as right $\O$-module on the projective line over $\Z$.

\textbf{Notation}. It follows from Definition \ref{principalparts} and Proposition \ref{local} that the 
sheaf of jets $\Pr^k(\O(n))|_{U_0}$ as a right $\O_{U_0}$-module equals
\[ \Z[u]\{1\otimes x_0^n,du\otimes x_0^n, \dots , du^k\otimes x_0^n\} \]
as a free $\Z[u]$-module, where $du^p=(t-u)^p$. 
It follows that  $\Pr^k(\O(n))|_{U_1}$ as a right $\O_{U_1}$-module equals
\[ \Z[v]\{1\otimes x_1^n,dv\otimes x_1^n, \dots , dv^k\otimes x_1^n\} \]
as a free $\Z[v]$-module, where $v=\frac{1}{u}$ and $dv^p=(\t-\u)^p$. 
Let 
\[ D=\{1\otimes x_0^n,du\otimes x_0^n, \dots , du^k\otimes x_0^n\} \]
and 
\[ D'=\{1\otimes x_1^n,dv\otimes x_1^n, \dots , dv^k\otimes x_1^n\} .\]
We wish to calculate the transition-matrix 
\[ R^k_n=[I]^{D'}_D\]
describing the right module structure on the sheaf of jets. 
The following theorem gives the structure-matrix $R^k_n$:
\begin{theorem} \label{rightstructure} Consider $\Pr^k(\O(n))$ as right $\O$-module on $\P$. If $n\in \Z$ and $k\geq 1$ 
the following formulas hold: 
\begin{equation} \label{right1}
 1\otimes x_1^n=u^n\otimes x_0^n 
\end{equation} 
and
\begin{equation}\label{right2}
 dv^p\otimes x_1^n=\sum_{j=p}^{k}(-1)^{j}\binom{j-1}{p-1}u^{n-p-j}du^{j}\otimes x_0^n, 
\end{equation}
where $1\leq p \leq k$. 
\end{theorem} 
\begin{proof} See \cite{maa2}.
\end{proof} 
%
%
%
%
%
The Theorems \ref{leftstructure} and \ref{rightstructure} describe completely the matrices $L^k_n$ and $R^k_n$ 
for all $n\in \Z$ and $k\geq 1$. By \cite{GRO} Theorem 2.1 and \cite{harder} Theorem 3.1 we know that any finite rank locally free
sheaf on $\P^1$ over any field splits into a direct sum of invertible sheaves. The formation of jets
commutes with direct sums, hence it follows from Theorems \ref{leftstructure} and \ref{rightstructure} that we have
calculated the transition matrix of $\Pr^k(\E)$ for any locally free finite rank sheaf $\E$ as left and right
$\O$-module on the projective line over any field. 

In the paper \cite{maa2} the structure of the sheaf of jets is studied and
the structure matrices $L^k_n,R^k_n$ are diagonalized explicitly. The
following structure theorem is obtained:

\begin{theorem} \label{theorem1} Let $k\geq 1$, and consider $\Pr^k(\O(d))$ as left $\O_\P$-module. 
If $k \leq d$ or $d<0$, there is an isomorphism
\begin{equation} \label{iso1} 
\Pr^k(\O(d))^L\cong \oplus^{k+1}\O(d-k)
\end{equation} 
of $\O$-modules. 
If $0 \leq d <k$ there is an isomorphism
\begin{equation}
\Pr^k(\O(d))^L\cong \O^{d+1}\oplus\O(d-k-1)^k
\end{equation}
as left $\O$-modules.
Let $b\in \mathbf{Z}$ and $k\geq 1$. There is an isomorphism
\begin{equation}
\Pr^k(\O(d))^R\cong \O(d)\oplus \O(d-k-1)^k
\end{equation}
as right $\O$-modules.
\end{theorem} 
\begin{proof}See \cite{maa2}.
\end{proof}

\section{On the left and right $P$-module structure on the projective line}

In this section we give a complete classification of the left and right
$P$-module structure of the jets on the projective line over a field of characteristic zero 
using the same techniques as in \cite{maa1}. 

 Let in general $V$ be a  finite dimensional vector space over a field $F$ of characteristic 
zero. Any affine algebraic group $G$ is a closed subgroup of $\GL (V)$ for some $V$, and given any closed subgroup
$H\subseteq G$, there exists a quotient map $G\rightarrow G/H$ with nice properties. The variety $G/H$ is smooth
and quasi projective, and the quotient map is universal with respect to $H$-invariant morphisms of varieties.
$F$-rational points of the quotient $G/H$ correspond to orbits of $H$ in $G$ (see \cite{jantzen}, section I.5).
Moreover: any finite-dimensional
$H$-module $\rho$ gives rise to a finite rank $G$-homogeneous vector bundle $E=E(\rho )$ and by \cite{AKH}, chapter 4 this correspondence
sets up an equivalence of categories between the category of linear finite dimensional representations of $H$ and 
the category of finite rank $G$-homogeneous vector bundles on $G/H$. There exists an equivalence of categories 
between the category of finite rank $G$-homogeneous vector bundles and the category of finite rank locally free 
sheaves with a $G$-linearization, hence we will use these two notions interchangeably. 

 Fix a line $L$ in $V$, and let $P$ be the closed subgroup of $\SL(V)$ stabilizing $L$. The quotient 
$\SL(V)/P$ is naturally isomorphic to projective space $\P$ parametrizing lines in $V$, and if we choose a basis $e_0,\cdots ,e_n$ 
for $V$, the quotient map
\[ \pi: \SL(V) \rightarrow \P \]
can be chosen to be defined as follows: map any matrix $A$ to its first column-vector. It follows that $\pi$ 
is locally trivial in the Zariski topology, in fact it trivializes over the basic open subsets $U_i$ of projective
space. One also checks that any $\SL(V)$-homogeneous vectorbundle on $\P$ trivializes over the basic open
subsets $U_i$. Any line bundle $\O(d)$ on projective  space is $\SL(V)$-homogeneous coming from a unique character 
of $P$ since $\SL(V)$ has no characters. 
The sheaf of jets $\Pr^k(\O(d))$ is a sheaf of bi-modules, locally free as left and right
$\O$-module separately. It has a left and right $\SL(V)$-linearization, hence we may classify the $P$-module
structure of $\Pr^k(\O(d))$ corresponding to the left and right structure, and that is the aim of this section. 

 The  calculation of the representations corresponding to $\Pr^k(\O(d))$ as left and right $P$-module is contained in Theorems
\ref{theorem1} and \ref{theorem2}. As a byproduct we obtain the
classification of the splitting type of the jet bundles obtained in
\cite{maa2}: this is Corollaries \ref{corollary1} and \ref{corollary2}.

Let in the following section  $V$ be a vector space over $F$ of dimension two. 
$\P=\SL(V)/P$ is the projective line parametrizing lines in $V$. There exists two exact sequences of $P$-modules.
\[ 0 \rightarrow L \rightarrow V \rightarrow Q \rightarrow 0 \]
and 
\[ 0\rightarrow m \rightarrow V^* \rightarrow L^* \rightarrow 0 ,\]
and one easily sees that $m\cong L$ as $P$-module. 

Let $p,q:\P\times \P\rightarrow \P$ be the canonical 
projection maps, and let $I\subseteq \O_{\P\times \P}$ be the ideal of the diagonal. Let $\O_{\Delta ^k}=\O_{\P\times \P}/I^{k+1}$
be the \emph{$k'th$ order infinitesimal neighborhood of the diagonal}.
Recall the definition of the sheaf of jets:

\begin{definition} Let $\E$ be an $\O_\P$-module. Let $k\geq 1$ and let 
\[ \Pr^k\E)=p_*(\O_{\Delta^k}\otimes q^*\E) \]
be the \emph{k'th order sheaf of jets of } $\E$. 
\end{definition}
If $\E$ is a sheaf with an $\SL(V)$-linearization, it follows
$\Pr^k(\E)$ has a canonical $\SL(V)$-linearization.
There is an exact sequence:
\begin{equation} \label{exact1} 
 0\rightarrow I^{k+1} \rightarrow \O_{\P\times \P} \rightarrow \O_{\Delta ^k} \rightarrow 0 .
\end{equation} 
Apply the functor $\R p_*(-\otimes q^*\O(d))$ to the sequence \ref{exact1} to obtain a long exact sequence of $\SL(V)$-linearized 
sheaves

\begin{equation} \label{exact2} 
0 \rightarrow p_*(I^{k+1} \otimes q^*\O(d)) \rightarrow p_*q^*\O(d) \rightarrow \Pr^k(\O(d))^L \rightarrow 
\end{equation} 
\[ \R^1 p_*(I^{k+1}\otimes q^*\O(d)) \rightarrow \R ^1p_*q^*\O(d) \rightarrow \R ^1p_*(\O_{\Delta ^k}\otimes q^*\O(d)) 
\rightarrow  \cdots \]
of $\O_\P$-modules. We write $\Pr^k(\O(d))^L$ to indicate that we use the left structure of the jets.
We write $\Pr^k(\O(d))^R$ to indicate right structure. 
The terms in the sequence \ref{exact2} are locally free since they are coherent and any coherent 
sheaf with an $\SL(V)$-linearization is locally free. 

\begin{proposition}\label{vanishing} Let $\E$ be an $\SL(V)$-linearized sheaf with support in $\Delta\subseteq \P\times \P$.
 For all $i\geq 1$ the following holds:
\[ \R ^ip_*(\E) =\R^iq_*(\E)=0. \]
\end{proposition} 
\begin{proof} Let $e\in \P$ be the distinguished point and consider the fiber diagram
\[
\diagram  \P \rto^j \dto^{\tilde{p},\tilde{q}} & \P\times \P \dto^{p,q} \\
           \Spec(\kappa(e)) \rto^i & \P \enddiagram \]
Since $\R^ip_*(\E)$ and $\R^iq_*(\E)$ have an $\SL(V)$-linearization it is enough to check the statement of the lemma on
the fiber at $e$. We get by \cite{HH}, Proposition III.9.3 isomorphisms 
\[ \R^ip_*(\E)(e)\cong \R^i\tilde{p}_*(j^*\E) \]
and 
\[ \R^iq_*(\E)(e)\cong \R^i\tilde{q}_*(j^*\E) ,\]
and since $j^*\E$ is supported on a zero-dimensional scheme, the lemma follows.
\end{proof} 
It follows from the Lemma that $\R ^1p_*(\O_{\Delta ^k}\otimes q^*\O(d))=\R ^1q_*(\O_{\Delta ^k}\otimes q^*\O(d))=0$ 
since $\O_{\Delta ^k}\otimes q^*\O(d)$ is supported on the diagonal. 
Hence we get an exact sequence of $P$-modules
when we pass to the fiber of \ref{exact2} 
at the distinguished point $e\in \P$. Let $\m\subseteq \O_\P$ be the sheaf of ideals of $e$. 
By \cite{HH} Theorem III.12.9 and Lemma \ref{vanishing} we get the following exact sequence
of $P$-modules:
\begin{equation} \label{exact3} 
0\rightarrow \H^0(\P, \m^{k+1}\O(d)) \rightarrow \H^0(\P,\O(d)) \rightarrow  \Pr^k(\O(d))(e) \rightarrow 
\end{equation}
\[ \H^1(\P, \m^{k+1}\O(d)) \rightarrow \H^1(\P, \O(d)) \rightarrow 0 .\]

\begin{proposition}\label{classification} Let $k\geq 1$ and $d < k$. Then there is an isomorphism
\[ \H^1(\P, \m^{k+1}\O(d)) \cong \sym^{k+1}(L)\otimes \sym^{k-d-1}(V) \]
of $P$-modules.
\end{proposition} 
\begin{proof} There is an isomorphism of sheaves $\O(-k-1)\cong \m^{k+1}$ defined as follows:
\[ x_0^{-k-1}\rightarrow t^{k+1} \] 
on the open set $D(x_0)$ where $t=x_1/x_0$. 
On the open set $D(x_1)$ it is defined as follows:
\[ x_1^{-k-1}\rightarrow 1 .\]
We get an isomorphism $\O(d-k-1)\cong \m^{k+1}\O(d)$ of sheaves, but the corresponding inclusion of
sheaves
\[ \O(d-k-1)\rightarrow \O(d) \]
is not a map of $P$-linearized sheaves, since it is zero on the fiber at $e$. Hence we must twist by the character
of $\m^{k+1}=\O(-k-1)$ when we use Serre-duality. We get
\[ \H^1(\P,\m^{k+1}\O(d))\cong \sym^{k+1}(L)\otimes \H^0(\P,\O(k-d-1))^*\cong  \]
\[ \sym^{k+1}(L)\otimes \sym^{k-d-1}(V), \]
and the proposition follows. 
\end{proof}

We first give a complete classification of the left $P$-module structure of the jets. The result is the
following.

\begin{theorem} \label{theorem1} Let $k\geq 1$, and consider $\Pr^k(\O(d))$ as left $\O_\P$-module. 
If $k \leq d$, there exists an isomorphism
\begin{equation} \label{iso1} 
\Pr^k(\O(d))^L(e)\cong \sym^{k-d}(L^*)\otimes \sym^k(V^*)
\end{equation} 
of $P$-modules. If $d\geq 0 $ there exist an isomorphism
\begin{equation} \label{iso2} 
\Pr^k(\O(-d))^L(e)\cong \sym^{k+d}(L)\otimes \sym^k(V)
\end{equation}
of $P$-modules. If $0\leq d <k$ there exist an isomorphism
\begin{equation} \label{iso3} 
\Pr^k(\O(d))^L(e)\cong \sym^d(V^*)\oplus \sym^{k+1}(L)\otimes \sym^{k-d-1}(V)
\end{equation}
of $P$-modules. 
\end{theorem} 
\begin{proof} The isomorphism from \ref{iso1} follows from Theorem 2.4 in \cite{maa1}. 
We prove the isomorphism \ref{iso2}: Since $-d<0$ we get an exact sequence
\[ 0\rightarrow \Pr^k(\O(-d))(e)\rightarrow \H^1(\P, \m^{k+1}\O(-d))\rightarrow \H^1(\P, \O(-d)) \rightarrow 0 \]
of $P$-modules. 
By proposition \ref{classification} we get the exact sequence of $P$-modules
\[ 0 \rightarrow \Pr^k(\O(-d))^L(e)\rightarrow \sym^{k+1}(L)\otimes \sym^{d+k-1}(V)\rightarrow \sym^{d-2}(V) \rightarrow 0 .\]
The map on the right is described as follows: dualize to get the following:
\[ \sym^{d-2}(V^*)\cong \sym^{k+1}(L^*)\otimes \sym^{k+1}(L)\otimes \sym^{d-2}(V^*)\cong  \]
\[  \sym^{k+1}(L^*)\otimes \sym^{k+1}(\m)\otimes \sym^{d-2}(V^*) \rightarrow \sym^{k+1}(L^*)\otimes \sym^{d+k-1}(V^*).\]
The map is described explicitly as follows:
\[ f(x_0,x_1)\rightarrow  x_0^{k+1}\otimes x_1^{k+1}f(x_0,x_1). \]
If we dualize this map we get a map
\[ \sym^{k+1}(L)\otimes \sym^{d+k-1}(V)\rightarrow \sym^{d-2}(V) ,\]
given explicitly as follows:
\[ e_0^{k+1}\otimes f(e_0,e_1)\rightarrow \phi(f(e_0,e_1)) ,\]
where $\phi$ is $k+1$ times partial derivative with respect to the $e_1$-variable. 
There exists a natural map
\[ \sym^{k+d}(L)\otimes \sym^k(V)\rightarrow \sym^{k+1}(L)\otimes \sym^{d+k-1}(V) \]
given by 
\[ e_0^{k+d}\otimes f(e_0,e_1)\rightarrow e_0^{k+1}\otimes e_0^{d-1}f(e_0,e_1) \]
and one checks that this gives an exact sequence
\[ 0\rightarrow \sym^{k+d}(L)\otimes \sym^{k}(V)\rightarrow \sym^{k+1}(L)\otimes \sym^{d+k-1}(V)\rightarrow \]
\[ \sym^{d-2}(V) \rightarrow 0 ,\]
hence we get an isomorphism 
\[ \Pr^k(\O(-d))^L(e)\cong \sym^{k+d}(L)\otimes \sym^k(V), \]
and the isomorphism from \ref{iso2} is proved.
We next prove isomorphism \ref{iso3}: By vanishing of cohomology on the projective line we get the following
exact sequence:
\[ 0\rightarrow \H^0(\P, \O(d)) \rightarrow \Pr^k(\O(d))^L(e) \rightarrow \H^1(\P, \m^{k+1}\O(d)) \rightarrow 0 .\]
We get by Proposition \ref{classification} an exact sequence of $P$-modules
\[0 \rightarrow \sym^d(V^*)\rightarrow \Pr^k(\O(d))^L(e) \rightarrow   \sym^{k+1}(L)\otimes \sym^{k-d-1}(V) \rightarrow 0 .\]
It splits because of the following:
\[ \Ext^1_P(\sym^{k+1}(L)\otimes \sym^{k-d-1}(V), \sym^d(V^*) )= \]
\[  \Ext^1_P(\rho, \sym^{k+1}(L^*)\otimes \sym^{k-d-1}(V^*)\otimes\sym^{k-1}(V^*)). \]
Here $\rho$ is the trivial character of $P$. 
Since there is an equivalence of categories between the category of $P$-modules and $\SL(V)$-linearized sheaves we get
again by equivariant Serre-duality 
\[\Ext^1_P(\rho, \sym^{k+1}(L^*)\otimes \sym^{k-d-1}(V^*)\otimes \sym^{k-1}(V^*))= \H^1(\P, \O(k+1)\otimes \E_1) =\] 
\[  \oplus^{e_1}\H^1(\P, \O(k+1))=\oplus^{e_1}\H^0(\P, \O(-k-3))^*=0 .\]
Here $\E_1$ is the $\SL(V)$-linearized sheaf corresponding to $\sym^{k-d-1}(V^*)\otimes \sym^{k-1}(V^*)$ and $e_1$ is the rank of
$\E_1$. 
Hence we get
\[ \Pr^k(\O(d))^L(e)\cong \sym^d(V^*)\oplus \sym^{k+1}(L)\otimes \sym^{k-d-1}(V), \]
and the isomorphism from \ref{iso3} is proved hence the theorem follows. 
\end{proof} 

As a corollary we get a result on the splitting type of the jets as left module on the projective line.

\begin{corollary} \label{corollary1} The splitting type of $\Pr^k(\O(d))$ as left $\O_\P$-module  is as follows:
If $k\geq 1$ and $d<0$ or $d\geq k$ there exist an isomorphism
\begin{equation}
\Pr^k(\O(d))^L\cong \oplus^{k+1} \O(d-k) 
\end{equation} 
of left $\O$-modules.
If $0\leq d <k$ there exist an isomorphism
\begin{equation}
\Pr^k(\O(d))^L\cong \O^{d+1}\oplus \O(-k-1)^{k-d} 
\end{equation} 
of left $\O$-modules.
\end{corollary}
\begin{proof} This follows directly from Theorem \ref{theorem1}. 
\end{proof} 

We next give a complete classification of the right $P$-module structure of the jets. The result is the following.

\begin{theorem} \label{theorem2}  Let $k\geq 1$ and consider $\Pr^k(\O(d))$ as right module. 
If $d>0$ there exist an isomorphism
\begin{equation} \label{iiso1}
\Pr^k(\O(-d))^R(e)\cong \sym^d(L)\oplus \sym^{k+d+1}(L)\otimes \sym^{k-1}(V)
\end{equation}
of $P$-modules.
If $d\geq 0$ there exist an isomorphism
\begin{equation} \label{iiso2}
\Pr^k(\O(d))^R(e)\cong \sym^{d}(L^*)\oplus \sym^{d-k-1}(L^*)\otimes \sym^{k-1}(V)
\end{equation} 
as $P$-modules.
\end{theorem}
\begin{proof} 
We prove the isomorphism \ref{iiso1}: Using the functor $\R^iq_*(-\otimes q^*\O(-d))$ we get a long exact sequence
of $\SL(V)$-linearized sheaves
\[ 0 \rightarrow q_*(I^{k+1})\otimes \O(-d) \rightarrow q_*q^*\O(-d) \rightarrow \Pr^k(\O(-d))^R \rightarrow \R^1q_*(I^{k+1})\otimes \O(-d)
\rightarrow \]
\[ \R^1q_*(\O_{\P\times \P})\otimes \O(-d) \rightarrow 0 .\]
It is exact on the right because of Lemma \ref{vanishing}. 
Here we write $\Pr^k(\O(d))^R$ to indicate we use the right structure of the jets. 
We take the fiber at $e\in \P$ and using Cech-calculations for coherent sheaves on the projective line, we obtain the following
exact sequence of $P$-modules:
\[ 0\rightarrow \H^0(\P, \O_\P)\otimes \O(-d)(e) \rightarrow \Pr^k(\O(-d))^R(e) \rightarrow \H^1(\P, \m^{k+1})\otimes 
\O(-d)(e) \rightarrow 0 \]
Hence we get by Proposition \ref{classification} the following exact sequence of $P$-modules:
\[ 0\rightarrow  \sym^d(L)\rightarrow \Pr^k(\O(-d))^R(e) \rightarrow \sym^d(L)\otimes \sym^{k+1}(L)\otimes \sym^{k-1}(V)\rightarrow 0 \]
It is split exact because of the following argument using $\Ext$'s and equivariant Serre-duality: 
\[ \Ext^1_P(\sym^{d+k+1}(L)\otimes \sym^{k-1}(V), \sym^d(L))=\Ext^1_P(\rho, \sym^{k+1}(L^*)\otimes \sym^k(V^*)),\]
where $\rho$ is the trivial character of $P$. By equivariant Serre duality we get 
\[\Ext^1_P(\rho, \sym^{k+1}(L^*)\otimes \sym^k(V^*))= \H^1(\P, \O(k+1)\otimes \E_2) = \] 
\[ =\oplus^{e_2} \H^0(\P, \O(-k-3))^* =0 .\]
Here $\E_2$ is the vectorbundle corresponding to $\sym^k(V^*)$ and $e_2$ is the rank of $\E_2$. 
Hence we get the desired isomorphism
\[ \Pr^k(\O(-d))^R(e)\cong \sym^d(L)\oplus \sym^{d+k+1}(L)\otimes \sym^{k-1}(V), \]
and isomorphism \ref{iiso1} is proved. 

We next prove the isomorphism \ref{iiso2}: Using the functor $\R^iq_*(-\otimes q^*\O(d))$ we get a long exact sequence
of $\SL(V)$-linearized sheaves
\[ 0 \rightarrow q_*(I^{k+1})\otimes \O(d) \rightarrow q_*q^*\O(d) \rightarrow \Pr^k(\O(d))^R \rightarrow \R^1q_*(I^{k+1})\otimes \O(d)
\rightarrow \]
\[ \R^1q_*(\O_{\P\times \P})\otimes \O(d) \rightarrow 0 .\]
We take the fiber at $e\in \P$ and using Cech-calculations for coherent sheaves on the projective line, we obtain the following
exact sequence of $P$-modules:
\[ 0\rightarrow \H^0(\P, \O_\P)\otimes \O(d)(e) \rightarrow \Pr^k(\O(d))^R(e) \rightarrow \H^1(\P, \m^{k+1})\otimes \O(d)(e) \rightarrow 0. \]
Proposition \ref{classification} gives the following sequence of $P$-modules
\[ 0 \rightarrow \sym^d(L^*) \rightarrow  \Pr^k(\O(d))^R(e) \rightarrow \sym^d(L^*)\otimes \sym^{k+1}(L)\otimes \sym^{k-1}(V)\rightarrow 0 \]
It splits because of the following $\Ext$ and equivariant Serre-duality argument:
\[ \Ext^1_P(\sym^d(L^*)\otimes \sym^{k+1}(L)\otimes \sym^{k-1}(V), \sym^d(L^*) )= \]
\[ \Ext^1_P(\rho, \sym^{k+1}(L^*)\otimes \sym^{k-1}(V^*)).\]
Let $\E_3$ be the vectorbundle corresponding to $\sym^{k-1}(V^*)$ and let $e_3$ be its rank. We get
\[ \Ext^1_P(\rho, \sym^{k+1}(L^*)\otimes \sym^{k-1}(V^*))= \H^1(\P, \O(k+1)\otimes \E_3 )= \]
\[ =\oplus^{e_3}\H^0(\P, \O(-k-3))^*=0, \]
hence we get the isomorphism
\[ \Pr^k(\O(d))^R(e)\cong \sym^d(L^*)\oplus \sym^{d-k-1}(L^*)\otimes \sym^{k-1}(V), \]
and isomorphism \ref{iiso2} follows. 
\end{proof}

As a corollary we get a result on the splitting type of the jets as right module.

\begin{corollary} \label{corollary2} Let $k\geq 1$ and $d\in \Z$. There exist an isomorphism
\begin{equation} \label{isoo}
\Pr^k(\O(d))^R\cong \O(d)\oplus \O(d-k-1)^k 
\end{equation}
of right $\O_\P$-modules.
\end{corollary}
\begin{proof} This follows directly from Theorem \ref{theorem2}. 
\end{proof}

Note that Corollary \ref{corollary1} and \ref{corollary2} recover Theorem \ref{theorem1}, hence we have used elementary
properties of representations of $\SL(V)$ to classify sheaves of 
bi-modules on the projective line over any field of characteristic zero. 

On projective space of higher dimension there is the following result:
Let $\p^N=\SL(V)/P$ where $P$ is the subgroup fixing a line, and let
$\O(d)$ be the linebundle with $d\in \mathbf{Z}$. It follows $\O(d)$
has a canonical $\SL(V)$-linearization.

\begin{theorem} \label{maintheorem} For all $1\leq k <d$, 
the representation corresponding to $\Pr^k(\O(d))^L$ is 
$\sym^{d-k}(L^*)\otimes \sym^k(V^*)$. 
\end{theorem}
\begin{proof} See \cite{maa1}.
\end{proof}

Note that the result in Theorem \ref{maintheorem} is true over any
field $F$ if $char(F)>n$.

\begin{corollary} For all $1\leq k<d$, $\Pr^k(\O(d))$ decompose as
  left $\O$-module as $\oplus^{\binom{N+k}{N} }\O(d-k)$.
\end{corollary}
\begin{proof} See \cite{maa1}.
\end{proof}

\textbf{Acknowledgements.} Thanks to Michel Brion and Dan Laksov for comments.


\noindent
      Helge Maakestad\\
      Dept. of Mathematics\\
      7491 Trondheim\\
      Norway\\
      maakesta@math.ntnu.no\\

\label{lastpage}
\end{document}